\documentclass{svjour3}

\usepackage{amscd,amsfonts,amssymb,amsmath,latexsym,array,hhline,xcolor,graphicx}

\usepackage{latexsym}
\usepackage{times}

\newcommand\F{\mbox{I\kern-2pt F}}

\newcommand\cF{{\cal F}}

\newcommand\e{{\varepsilon}}
\newcommand{\wt}{\widetilde}

\def\E{{\bf E}}
\def\P{{\bf P}}
\def\Chi{{\bf 1}}

\def\bbr{{\mathbb R}}

\newcommand\fdem{$\Box$}
\newcommand\beq{\begin{equation}}
\newcommand\eeq{\end{equation}}
\newcommand\bea{\begin{eqnarray}}
\newcommand\eea{\end{eqnarray}}
\newcommand\bean{\begin{eqnarray*}}
\newcommand\eean{\end{eqnarray*}}

\begin{document}
\title{In the Life Insurance Business Risky Investments are Dangerous}

\author{Yuri Kabanov  
\and Serguei\ ~ Pergamenshchikov}

\date{}
\institute{
  \at
              Laboratoire de Math\'ematiques, Universit\'e de
Franche-Comt\'e, \\ 16 Route de Gray,
 25030 Besan\c{c}on, cedex,
France, and \\
International Laboratory of Quantitative Finance, National Research University  Higher School of Economics,\\   Moscow, Russia\\
 \email{Youri.Kabanov@univ-fcomte.fr}    \\
\and Laboratoire de Math\'ematiques Rapha\"el Salem, \\
Universit\'e de Rouen, 
Avenue de l'Universit\'e, BP.12
Technop\^ole du Madrillet
F76801 Saint-\'Etienne-du-Rouvray, France,  and \\
International Laboratory of Quantitative Finance, National Research University  Higher School of Economics,\\
   Moscow, Russia\\
  \email{Serge.Pergamenchtchikov@univ-rouen.fr}   }

\date{Received: date / Accepted: date}

\titlerunning{In the Life Insurance Business Risky Investments are Dangerous}

\maketitle

\qquad \qquad \qquad \qquad \qquad \qquad \qquad \qquad \qquad \quad
{\sl Dedicated to the memory of Marc Yor.}\\
\\

\begin{abstract}
We investigate  models of the  life annuity insurance when the company invests its reserve into a risky asset with  price following a geometric
Brownian motion.  Our main result is an exact asymptotic of the ruin probabilities for  the case of exponentially distributed benefits. As in the case of non-life insurance with exponential claims, the ruin probabilities are either decreasing with a rate given by a power function (the case of small volatility)
or equal to unit identically (the case of large volatility). The result allows us 
to quantify the share of reserve  to invest into such a risky asset to avoid a catastrophic outcome: the ruin with probability one. We address also the question of smoothness of the ruin probabilities as a function of the initial reserve for generally  distributed jumps. 
\end{abstract}




 \section{Introduction} 
 In the modern world insurance companies operate in a financial environment. In particular, they 
 invest their reserves into various assets and this may add more risk to the business. 
 To our knowledge, the first  model of an insurance company investing its capital into a risky asset 
 has appeared in the short note \cite{Fr} where the author provided arguments showing that the asymptotic behavior of ruin probability is radically different from that in the  
classical Lundberg--Cram\'er model.   A rigorous analysis in \cite{FrKP} confirmed the conjecture. 
Using Kalashnikov's estimates for linear finite difference equations with random coefficients, \cite{Kalash},  it was shown that, independently of the safety loading,  the ruin is imminent with unit probability when  the volatility $\sigma$ of the stock price is large with respect to the instantaneous rate of return $a$ (namely, when $2a/\sigma^2<1$), and the ruin probability is decreasing as a power function, when the volatility is small (namely, when $2a/\sigma^2>1$). For the model with exponentially distributed claims the exact asymptotic was found.  
Threshold case $2a/\sigma^2=1$ was studied in the papers \cite{PZ}, \cite{PerErr} where it was shown, using techniques based on a renewal theorem, that the ruin is imminent with unit probability.  
The setting of  \cite{FrKP} and \cite{PZ} is of the so-called non-life insurance: the company receives a flow of contributions and pay primes. The ruin occurs when a new claim arrives and its value  
is too large to be covered by the reserve: the risk process exits from the positive half-axis by a jump. The model can be studied in the discrete-time framework but the method of differential equations happens to be more efficient in the case of exponential claims where it  allows to get the exact asymptotic.     

In the present note we investigate the setting of life annuity, or pension,  insurance, when the company pays to the policyholder a rent and 
earns a premium when the insured person dies, see, e.g., \cite{Gr}, p. 8 or \cite{Sax}; in the classical literature  such models are called {\it models with negative risk sums}.    
The highly stylized model of the collective risk theory may describe the situation where annuitant, entering the contract, pays a lump sum, e.g., the savings during his working life, getting in exchange regular payments until the death. One can interpret the premium as the reserve release. If the company has a large "stationary" portfolio of such contracts, with new arriving customers,  then it is reasonable to think that the dates of reserve releases form a Poisson process and the sizes of premiums are independent random variables of the same distribution. A special type of life annuity, called {\it le viager} exists in France. This practice is, basically, an exchange of a property, a house or an apartment, for life-long regularly payments. The 
annuitant enjoys a regularly income without having to give a house. For the company having a portfolio of viagers  the amount of premiums depends on the real estate market.

For this classical model, where the reserve process has jumps upwards and may leave the positive half-axis only in a continuous way,  the ruin problem  can be easily reduced to the ruin problem for the non-life insurance using the so-called duality method, \cite{Asm}.  Its idea is to define the "dual" process by replacing the line segments of the original process between consecutive jumps by downward jumps of the same depth and the upward jumps by  line segments   of the same height with positive fixed slope. Note that in the  literature the models with upward jumps are often referred to as {\it dual models}, \cite{AGS}, \cite{ABL}. 

The duality method does not work in our setting where the capital of the company or its fraction is invested into a risky asset. The change of two signs to the opposite ones in the equation  defining the  dynamics of the reserve leads to technical complications. In particular, now the ruin  may happen before the instant of the first jump and the latter is no more  the instant of the regeneration, after which the process starts afresh provided being strictly positive.       

Nevertheless, a suitable modification of arguments of    \cite{FrKP} and \cite{PZ} combined with  some ideas   allows us  to obtain the  asymptotic of the ruin probability as the initial reserve tends to infinity. Expectedly,  it is  the same as in the non-life insurance case. In many countries there are rules allowing  insurance companies to invest only a small share of their reserves into risky assets. Our simple  model confirms that it is reasonable and even provides 
a quantitative answer.  To avoid the situation when the ruin happens with probability one,  the proportion of investment into the risky asset should be strictly less than $2a/\sigma^2$.   

It should be emphasized that the life insurance case is rather different and its study is not a straightforward exercise.  The main difficulty in deriving the integro-differential equation    
is to prove the smoothness of the ruin probability and the integrability of derivatives.  This issue  is already delicate in the non-life insurance case. 
Unfortunately,  it was not discussed it in  \cite{FrKP} where the reader was directed towards the literature.   Now we are not sure that we disposed at that time a clean reference.  The smoothness of the exit probability is discussed in many papers, see, e.g.,  an interesting article \cite{WangWu}
where it  the explicit formula for an exponential functional of Brownian motion due to Marc Yor is used. Unfortunately, the needed smoothness property was established only under constraints on coefficients. One of the authors of the present note  had a fruitful discussion with Marc on  a possibility to deduce smoothness of the ruin probability function and the integrability of its derivatives without using complicated explicit formulae. His suggestions are realized here for the life insurance case. Note that in the literature on non-life insurance one can find other methods  to  establish the smoothness, for example an approach similar to verification theorems in the stochastic control theory, see \cite{Belkina}.

The structure of the paper is as follows. Section 2 contains the formulation of the main results. In Section 3 we establish  an upper asymptotic bound for the exit  probability (from $]0,\infty[$)  for a solution of a non-homogeneous linear stochastic equation and a lower asymptotic bound for a small volatility case. As a tool we use the  Dufresne theorem rather than Goldie's renewal theorem from \cite{Go91}, which was the kew ingredient of arguments in \cite{PZ}.  
The proof of Theorem \ref{main.1}
asserting  that in the case of large volatility the ruin is imminent is given in Section 4. The regularity of the non-ruin probability $\Phi$ is studied in Section 5 using a method based on integral representations. At the end of this section we derive the integro-differential equation for $\Phi$. Section 6 contains the proof of the main theorem. Finally, in the appendix we provide a formulation of an ergodic theorem for an autoregression with random coefficients proven in \cite{PZ}.  
  
Kalashnikov's approach (developed further in his joint work with Ragnar Norberg \cite{Kalash-Nor}) plays an important role in our study. Our
techniques is elementary. More profound and general results can be
found in \cite{Paul-98}, \cite{Paul}, \cite{Paul-G}, \cite{Norberg}, \cite{Nyrh-99}, \cite{Nyrh-01} {\it et al}.

\section{The model}

We are given a stochastic basis $(\Omega,\cF,{\bf F}=(\cF_t)_{t\ge 0},\P)$ with a Wiener process $W$
independent of the  integer-valued random measure $p(dt,dx)$ with
the compensator $\tilde p(dt,dx)$.

Let us consider a process $X=X^u$ of the form
 \beq
 \label{risk}
 X_t=u+  a \int_0^t  X_s  ds+
  \sigma \int_0^t X_s  dW_s  -ct  +\int_0^t\int x
 p(ds,dx),
 \eeq
where  $a$ and $\sigma$ are arbitrary constants  and $c\ge 0$.

We shall assume that ${\tilde p}(dt,dx)=\alpha dtF(dx)$ where
$F(dx)$ is a probability distribution on $]0,\infty[$. In this
case the integral with respect to the jump measure is simply a
compound Poisson process. It can be written as
$\sum_{i=1}^{N_t}\xi_i$ where $N$ is a  Poisson process with
intensity $\alpha$ and $\xi_i$ are random variables with common
distribution $F$;  the random variables $W$, $N$, $\xi_i$, $i\in {\bf N}$, are
independent. We  denote by $T_n$ the successive jumps of $N$; the inter-arrival 
time intervals  
$T_i-T_{i-1}$ are   independent and exponentially distributed with the parameter $\alpha$. 

In our main result (Theorem \ref{main}) we assume that  $F$ is the
exponential distribution with parameter $\mu$.

 Let $\tau^u:=\inf \{t: X^u_t\le 0\}$ (the instant of ruin),
 $\Psi (u):=P(\tau^u<\infty)$ (the ruin probability),
 and $\Phi(u):=1-\Psi (u)$
 (the survival probability).

The parameter values $a=0$, $\sigma=0$, correspond to the
life insurance ( i.e. the dual) version of the 
Lundberg--Cram\'er model for which  the risk process is usually
written as 
\beq
\label{r}
r_t:=u  - ct  +\sum_{i=1}^{N_t}\xi_i.
\eeq
 In the considered case the capital evolves due to continuously outgoing cash
flow with rate $c$ and incoming  random payoffs $\xi_i$ at times
forming an independent  Poisson process $N$ with intensity
$\alpha$. For the classical model with positive safety loading and $F$
having a ``non-heavy" tail, the Lundberg inequality provides an
encouraging information: the ruin probability decreases
exponentially as the initial capital $u$ tends to infinity.
Moreover, for the exponentially distributed claims the ruin
probability admits an explicit expression, see \cite{Asm} or
\cite{Gr}.

The more realistic case $a>0$, $\sigma=0$, corresponding to
non-risky investments, does not pose any problem 
 (see, e.g., \cite{Harrison}  for estimates of the exit probabilities covering this case). 
 
 We study here the case  $\sigma>0$. Now the equation (\ref{risk})
describes the evolution of the  reserve of an insurance company 
which pays a rent and continuously reinvests its capital  into an asset with the price
following a geometric Brownian motion.  The same model can be used 
for the description of the capital of a venture company funding R\&D and selling innovations, \cite{Bayraktar-Egami}.   

\smallskip
\noindent
{\bf Notations.} Throughout the paper we shall use the following abbreviations: 
$$
\kappa:=a-\frac12 \sigma^2, \qquad \beta:= \frac{2\kappa}{\sigma^2}=\frac{2a}{\sigma^2}-1, \qquad \eta_t:=\kappa t+\sigma W_t.   
$$  

\smallskip

The solution of the linear stochastic  equation (\ref{risk}) can be written using the  Cauchy formula: 
 \beq
 X_t=e^{\eta_t} \left(u+\int_{[0,t]}e^{-\eta_s}dr_s\right). 
 \eeq

 \begin{theorem}
\label{main}
  Let $F(x)=1-e^{-x/\mu}$, $x>0$. Assume that $\sigma>0$.

 (i) If
  $\beta>0$, then for some $K>0$
 \beq
  \Psi (u)= Ku^{-\beta}(1+o(1)), \quad u\to \infty.
  \eeq

(ii) If $\beta\le 0$, then $\Psi (u)=1$ for all $u>0$.
 \end{theorem}

The formulation of this theorem is exactly the same as in \cite{FrKP}
for the non-life insurance model (the case $\beta=0$ was analyzed in 
\cite{PZ}, \cite{PerErr}). 

One must admit that in the life insurance models, especially, in the case of a company keeping a portfolio of viager contracts, the hypothesis that the benefits (i.e. prices of houses) follow the exponential distribution, are highly unrealistic. Nevertheless, we can claim,  without any assumption on the distribution, that the ruin probabilities lay between two power functions, see Propositions \ref{main-0} and \ref{Pr.sec:Lb.1}.   

The next result (implying the statement $(ii)$ above) says 
that for $\delta>0$ the ruin is imminent. It requires only that the distribution $F$ has moments of positive order. 

\begin{theorem}
\label{main.1}
  Assume that there is $\delta>0$ such that 
$\E\xi^{\delta}_{1}<\infty$.  If $\beta\le 0$,   then $\Psi (u)=1$ for any  $u>0$.
 \end{theorem}

 The same model serves well in the situation where  only a fixed
part $\gamma\in ]0,1]$ of the capital is invested in the risky
asset: one should only replace the parameters $a$ and $\sigma$ in
(\ref{risk}) by $a\gamma$ and $\sigma\gamma$.  
Theorem \ref{main} implies that the ruin with probability one will be avoided only if $2a\gamma/(\sigma\gamma)^2>1$, i.e. when the share of investment into the risky asset is strictly less than $2a/\sigma^2$. 
\smallskip

It is worth  mentioning that our conclusion are robust and holds for more  general models. The reader may contest that intensity of outgoing payments $c$ is constant. 
Indeed, after the death of the annuitant the payments stop and the intensity must decrease while with a new customer it increases. Easy comparison arguments show that  the above statements hold in a more realistic situation where  $c=(c_t)$ is a random process such that $0<C_1\le c\le C_2$ where   $C_1$ and $C_2$ are constants.

The crucial part of the asymptotic analysis in Theorem 1 is based on the fact  
that for the Markov process given
by (\ref{risk}) the non-exit probability $\Phi (u)$ is smooth and  satisfies the
following equation:
\beq
\label {int-diff} \frac 12 \sigma^2u^2 \Phi''(u)+ (au-c)\Phi'(u) -
\alpha\Phi(u)+\alpha \int_0^\infty\Phi(u+y)dF(y)=0.
\eeq

With $\sigma>0$  this equation is of the second order and, hence,
requires two boundary conditions ---  in contrast to the classical case
($a=0$, $\sigma=0$) where it degenerates to an equation of the
first order requiring  a single boundary condition, see \cite{Gr}.  The estimate 
given in Proposition \ref{main-0} shows that $\Phi (\infty)=1$.  

There is an extensive literature concerning of the regularity of the survival probability for process with jumps  in the context of non-life insurance models and models based on the L\'evy processes, see, e.g., \cite{WangWu},  \cite{Ga-Gr}, 
 \cite{Granditz}, \cite{BelkinaKK}, \cite{BelKK}.  Our Theorem \ref{Th.sec:Reg.1}  requiring only the smoothness of $F$ and the integrability of $F''$ seems to be the first result on the regularity of the survival probability in the considered setting.
  
\section{Asymptotic bounds for the small volatility case}
\subsection{Upper asymptotic bound}

Let us consider the exit probability problem for a more general process $X=X^u$ of the form
 \begin{equation}
  \label{risk-0}
 X_t=u+  a \int_0^t  X_s  ds+
  \sigma \int_0^t X_{s}\,dW_{s}  - ct  +
  Z_{t},
 \end{equation}
where  $a, \sigma, c >0$ are arbitrary constants and 
$Z=(Z_{t})_{t\ge 0}$ is an increasing adapted c\`adl\`ag process starting from zero.  


 \begin{proposition}
\label{main-0}
Assume that  in the general model (\ref{risk}) the parameter $\beta>0$. 
Then 
$$
\limsup_{u\to\infty}\,u^{\beta}\,\Psi(u)\,\le\,\frac{2^{\beta}\,c^\beta}{\sigma^{2\beta}\beta\Gamma(\beta)}.
$$
 \end{proposition}
{\sl Proof.}  Let $Y$  be a solution of  the linear stochastic equation 
$$
Y_{t}=u+a\int^{t}_{0}\,Y_{s}\,d s+\sigma\,\int^{t}_{0}\,Y_{s}\,d W_{s}-c t.   
$$   
 Introducing the notation 
\beq
\label{defR}
R_t:=c\int^{t}_{0}\,e^{-\eta_v}\,d v 
\eeq
we express the solution as   
$$
Y_t:=e^{\eta_t}\big(u-R_t). 
$$ 
The  difference $X-Y$ satisfies  the linear equation with zero initial condition. Since   $Z$ is increasing,  we have the inequality $X\ge Y$ showing that the exit of $X$ from $]0,\infty[$ implies the exit of $Y$. Thus, 
\beq
\label{sec:Up.2-0}
\Psi(u)\le P(R_\infty>u). 
\eeq
and the asymptotic behavior of the ruin probability for this general model can be estimated by the tail behavior of the distribution function of  
$R_\infty$. Using the change of variable $v=(4/\sigma^2) t$ and observing that $B_t:=-(1/2)\sigma W_{(4/\sigma^2)t}$ is a Wiener process we obtain the representation  
$$
R_\infty=c\int_0^\infty  e^{-(a-\sigma^2/2)v-\sigma W_v}\,dv= 
\frac {4c}{\sigma^2} \int_0^\infty  e^{-2\beta t+2B_t}\,dt=:\frac {4c}{\sigma^2} A_\infty^{(-\beta)}. 
$$
The Dufresne theorem (see, \cite{Dufresne} or \cite{MatsumotoYor}, Theorem 6.2) claims that $A_\infty^{(-\beta)}$ is distributed as the random variable $1/(2\gamma)$ where $\gamma$ has the gamma distribution with parameter $\beta$. Thus, 
\bean
\P(R_\infty>u)&=&\P(2c/(\sigma^2\gamma) >u))=\P(\gamma<2c/(\sigma^2u))=
\frac{1}{\Gamma(\beta)}\int_0^{2c/(\sigma^2u)}x^{\beta-1}e^{-x}\,dx\\&\sim &
\frac{2^{\beta}\,c^\beta}{\sigma^{2\beta}\beta\Gamma(\beta)}u^{-\beta}
\eean
and the result follows from (\ref{sec:Up.2-0}). \fdem 

 \subsection{Lower asymptotic bound}
\label{sec:Lb}
The next result shows that the ruin probability decreases as the initial capital tends to infinity not faster than a certain power function.  

 \begin{proposition}
  \label{Pr.sec:Lb.1}
  Assume that   $\beta>0$.  
  Then there exists $\beta_{*}>0$ such that
 $$
 \liminf_{u\to\infty}\,u^{\beta_{*}}\,\Psi(u)\,>\,0\,.
 $$
 \end{proposition}
 \noindent {\sl Proof.}   
Let $Y=(Y_{k})_{k\ge 1}$ be the imbedded discrete-time process, that is the sequence of random variables defined recursively  as follows: 
\begin{equation}
\label{sec:Lb.1}
Y_{k}=M_{k} Y_{k-1}+Q_{k},\qquad Y_{0}=u,
\end{equation}
where 
$$
M_{k}=e^{\eta_{T_{k}}-\eta_{T_{k-1}}} 
\quad\mbox{and}\quad 
Q_{k}=\xi_{k}-c\,\int^{T_{k}}_{T_{k-1}}\, e^{\eta_{T_{k}}-\eta_{v}}\,d v\,.
$$
Let $\theta^u:=\inf\{k:\ Y_k\le 0\}$. 
It is clear that $X_{T_{k}}= Y_{k}$ for any $k\ge 1$. So, for any $u>0$
\begin{equation}
\label{sec:Lb.1-1}
\Psi(u)= \P(\tau^u<\infty) \ge 
\P(\theta^u<\infty).
\end{equation}
Take  $\varrho\in ]0,1[$ and chose   $B$ sufficiently large to ensure that
$$
B_{1}=B-\frac{1}{\varrho^{2}(1-\varrho)}>0. 
$$
Define the sets  
\begin{equation}\label{sec:Lb.2}
\Gamma_{k}=\{M_{k}\le \varrho\}\cap\{Q_{k}\le \varrho^{-1}\}\,,
\quad
D_{k}=\{M_{k}\le \varrho^{-1}\}\cap\{Q_{k}\le -B\}.
\end{equation}
\noindent On the set $\cap^{n}_{k=1}\,\Gamma_{k}$
we have
$$
Y_{n}=u\prod^{n}_{j=1}\,M_{j}
+\sum^{n}_{k=1}\,Q_{k}\,\prod^{n}_{j=k+1}\,M_{j}
\le u \varrho^{n}+\frac{1}{\varrho(1-\varrho)}.
$$
Therefore, on the set $\cap^{n}_{k=1}\,\Gamma_{k}\cap D_{n+1}$
$$
Y_{n+1}=M_{n+1} Y_{n} + Q_{n+1}\le u \varrho^{n-1}+\frac{1}{\varrho^{2}(1-\varrho)}
-B= u \varrho^{n-1} -B_1.
$$
It is easy to check that  $u \varrho^{n-1}\le  B_1$, when 
$u>B_1$ and  
$$ 
n=3+\left[\frac{\ln(u/B_{1})}{\vert \ln\varrho \vert} \right],
$$
where $[...]$ means the integer part.  
Therefore, 
$$
\P(\theta^u<\infty)\ge \P\left( \cap^{n}_{k=1}\,\Gamma_{k}\cap D_{n+1}\right)
=\left(\P(\Gamma_{1})\right)^{n}\P(D_{1}).
$$
Taking into account that 
$\P(\Gamma_{1})>0$ and 
$\P(D_{1})>0$, 
we obtain that
$$
\lim_{u\to\infty}\,u^{\beta_*} \P(\theta^u<\infty)=\infty
$$
for any 
$$
\beta_*>\frac{\ln\P(\Gamma_{1})}{\ln\varrho}.
$$
This implies the claim.
\fdem

\section{Large volatility: proof of Theorem \ref{main.1}}

We consider separately two cases: $\beta<0$ and $\beta=0$ and show that in both cases the ruin probability is equal to one. 

 \begin{proposition} 
 \label{Pr.sec:Lb.2}
  Assume that    $\beta<0$ and $\E\,\xi_1^\delta<\infty$ for some $\delta>0$.    Then  $\Psi (u)=1$ for all $u>0$. 
 \end{proposition}
 \noindent {\sl Proof.} As in the proof of Proposition \ref{Pr.sec:Lb.1} we consider the imbedded discrete-time process $Y=(Y_{k})_{k\ge 1}$ defined by (\ref{sec:Lb.1}). By virtue of (\ref{sec:Lb.1-1}) 
 it is sufficient to show that $\P(\theta^u<\infty)=1$. 
 Note that for $\delta \in ]0,-\beta[$
 $$
 \E e^{\delta \eta_t}= \E e^{\delta (\kappa t+\sigma W_t)}=e^{\delta t(\beta+\delta)\sigma^2/2}<1, 
 $$
 and, therefore,
 \beq
 \E M_{1} ^{\delta}=\int_0^\infty   \E e^{\delta \eta_t}\alpha e^{-\alpha t}dt <1. 
 \eeq
\noindent 
According to (\cite{BoSa}, p. 250)
 if a random variable $\nu$ is independent of $W$ and has the exponential distribution with parameter $\alpha$, then 
\begin{equation}
\label{sec:Up.4}
\E\,e^{q\max_{v\le \nu}(\mu v+W_{v})}
=\frac{\sqrt{2\alpha+\mu^{2}}-\mu}{\sqrt{2\alpha+\mu^{2}}-\mu-q}.
\end{equation}
provided that 
$$
\sqrt{2\alpha+\mu^{2}}-\mu-q>0. 
$$

Changing the variable and estimating the integrand by its maximal value we get that  
\bean
\E\Big(\int^{T_{1}}_{0}\, e^{\eta_{T_{1}}-\eta_{v}}\,d v \Big)^{\delta}
&=&
\E\Big(\int^{T_{1}}_{0}\, e^{\kappa v+\sigma W_{v}}\,d v \Big)^{\delta}
\le 
\E\, T^{\delta}_{1} e^{\delta \max_{v\le T_{1}}(\kappa v+\sigma W_{v})}\\
&=&\int _0^\infty t ^\delta\E\, e^{\delta \max_{v\le t}(\kappa v+\sigma W_{v})}\alpha e^{-\alpha t}d\,t\\
&\le&
2 \sup_{t\ge 0} \{t^{\delta}\,e^{-\alpha t/2}\} 
\E\, e^{\delta \max_{v\le \nu'}(\kappa v+\sigma W_{v})}\\
&=&
2e^{\delta(\ln(2\delta/\alpha)-1)}\,
\E\, e^{\delta \max_{v\le \nu'}(\kappa v+\sigma W_{v})},
\eean
where $\nu'$ is an  exponential random variable with parameter $\alpha/2$. In  view of the equality (\ref{sec:Up.4})
$$
\E \,e^{\delta \max_{v\le \nu'}(\kappa v+\sigma W_{v})}\,=\,
\frac{\sqrt{\sigma^{2}\alpha+\kappa^{2}}-\kappa}{\sqrt{\sigma^{2}\alpha+\kappa^{2}}-\kappa-\delta\sigma^{2}},
$$
provided that 
$$
\delta<\frac{\sqrt{\sigma^{2}\alpha+\kappa^{2}}-\kappa}{\sigma^{2}}, 
$$ 
i.e. for such $\delta>0$ we have the bound 
\beq
\label{sec:Up.4-0}
\E\Big(\int^{T_{1}}_{0}\, e^{\eta_{T_{1}}-\eta_{v}}\,d v \Big)^{\delta}\,
\le\,
2\,e^{\delta(\ln(2\delta/\alpha)-1)}\,
\frac{\sqrt{\sigma^{2}\alpha+\kappa^{2}}-\kappa}{\sqrt{\sigma^{2}\alpha+\kappa^{2}}-\kappa-\delta\sigma^{2}}\,.
\eeq

\noindent 
Using these estimates and the assumption of the proposition we conclude that  
$\E\, |Q_{1}|^{\delta}<\infty$ for sufficiently small $\delta>0$.  Thus, the hypothesis of the ergodic theorem for autoregression with random coefficients is fulfilled (see
Proposition  \ref{Pr.sec:A.1}). The latter claims that     
 for any bounded uniformly continuous function $f$
\begin{equation}
\label{sec:Up.5}
\P\hbox{-}\lim_{N}\frac 1N\sum^{N}_{k=1}\,f(Y_{k})\,=\,
\E\,f(\zeta),
\end{equation}
where 
$$
\zeta=Q_{1}\,+\sum^\infty_{k=2}\,Q_{k}\,
 \prod^{k-1}_{j=1}\, M_{j}.
 $$
Let us represent $\zeta$ in the form 
$$
\zeta=\xi_1-\int_0^{T_1}e^{\eta_{T_1}-\eta_v}dv+e^{\eta_{T_1}}\zeta_1,  \qquad \zeta_1:=\sum^\infty_{k=2}\,Q_{k}\,
 \prod^{k-1}_{j=2}\, M_{j}. 
$$ 
Clearly, the random variables $\xi$, $\zeta_1$ and $(\eta_{T_1},\int _0^{T_1}e^{\eta_{T_1}-\eta_v}dv)$ are independent. 
Moreover, Lemma \ref{support} given after the proof implies that the support of conditional distribution of the integral  $\int_0^{t}e^{\eta_{t}-\eta_v}dv$ given $\eta_t=x$ is unbounded from above. 
From this we easily infer that the support of distribution of $\zeta$ is unbounded from below. 
Thus, for 
$$
f(x)=\Chi_{\{x\le -1\}}+\vert x\vert\,\Chi_{\{-1<x<0\}}
$$
 the right-hand side of (\ref{sec:Up.5}) is strictly positive and, therefore,  
 $\P\big(\inf_{k\ge 1}\,Y_{k}<0\big)=1$. \fdem 

\begin{lemma}
\label{support} Let $\sigma>0$. Then  the support of conditional distribution of the random variable  
$$
I=\int_0^1 e^{\sigma W_v}d\,v
$$
given $W_1=y$ is unbounded from above. 
\end{lemma}
{\sl Proof.} It is well-known (see, e.g., \cite{RY}) that the conditional distribution of the Wiener process 
$(W_v)_{ v\le 1}$ given $W_1=y$ coincides with the distribution of the Brownian bridge $B^y$ with  $B^y_t =W_t+t(y-W_1)$. Thus, the conditional distribution of $I$ is the same as the unconditional distribution of  
$$
\tilde I:=\int_0^1 e^{\sigma (W_v+v(y-W_1))}d\,v.
$$
  Since 
the Wiener measure has the full support in the space $C_0[0,1]$ of continuous functions on $[0,1]$ with zero initial value, the support of distribution of $\tilde I$ is unbounded from above.  \fdem 




\begin{proposition} 
\label{Pr.sec:Lb.3}
   Assume that    $\beta=0$ and $\E\,\xi_1^\delta<\infty$ for some $\delta>0$.    Then  $\Psi (u)=1$ for all $u>0$. 
 \end{proposition}
 {\sl Proof.} In  the considered case the  imbedded discrete-time process is defined by  (\ref{sec:Lb.1}) with   
$$
M_{k}=e^{\sigma \Delta V_{k}} 
\quad\mbox{and}\quad 
Q_{k}=\xi_{k}-c \int^{T_{k}}_{T_{k-1}}\, e^{\sigma ( W_{T_{k}}-W_{v})}\,d v,
$$
where $V_k=W_{T_{k}}$ and $\Delta V_k:=V_k-V_{k-1}$.  
To study the asymptotic properties of the equation (\ref{sec:Lb.1})
we use the approach proposed in \cite{PZ} for the non-life insurance models. 
To this end, define recursively  the sequence of random variables putting $\theta_0:=0$ and  
\begin{equation}
\label{sec:Up.6}
\theta_{n}:=\inf\{k>\theta_{n-1}\colon\ V_k-V_{\theta_{n-1}} <0 \}, \qquad n\ge 1. 
\end{equation}
Note that    
$\theta_{n}=\sum^{n}_{j=1}\,\Delta \theta_j $ where $(\Delta \theta_j)_{j\ge 1}$ is a sequence 
of i.i.d. random variables distributed as $\theta_1$. 
It is known (see, e.g., XII.7, Theorem $1a$ in \cite{Fel}) that
\begin{equation}
\label{sec:Up.6-00}
C:=
\sup_{n\ge 1}\,n^{1/2}\P(\theta_{1}>n)\,<\,\infty
\,.
\end{equation}
Putting $y_{k}=Y_{\theta_{k}}$,
we obtain from  (\ref{sec:Lb.1}) that 
\begin{equation}
\label{sec:Up.7}
y_{k}={a}_{k}\,y_{k-1}+{b}_{k},\quad {y}_{0}=u,
\end{equation}
where
$$
{a}_{k}=\prod^{\Delta \theta_{k}}_{j=1} M_{\theta_{k-1}+j}
=
e^{\sigma (V_{\theta_k}-V_{\theta_{k-1}})}
$$
and
$$
{b}_{k}=
\sum^{\Delta \theta_k}_{l=1}
\,
\left( 
\prod^{\Delta \theta_k}_{j=l+1}\,M_{\theta_{k-1}+j}
\right)
\,Q_{\theta_{k-1}+l}.
$$
It is clear that ${a}_{k}<1$ a.s. Moreover, the first condition in Theorem~\ref{main.1}
and the inequality
(\ref{sec:Up.4-0}) with $\kappa=0$
implies that $\E \vert Q_{1} \vert^{\delta}<\infty$ for any sufficiently small $\delta$. Now, taking into account that 
$$
\vert {b}_{1} \vert \le \sum^{\Delta \theta_{1}}_{l=1}
\left( 
\prod^{\Delta \theta_1}_{j=l+1}M_{j}
\right)
\vert Q_{l}\vert
=
\sum^{\Delta \theta_{1}}_{l=1}
\,
\frac{{a}_{1}}{\prod^{l}_{j=1}\,M_{j}}
\,
\,\vert Q_{l}\vert
\,
\le\, 
\sum^{\Delta \theta_{1}}_{l=1}
\vert Q_{l}\vert ,
$$
we can get, for  $r\in ]0,1[$ and an increasing sequence of integers $l_n$, that 
\begin{align*}
\E\,\vert {b}_{1} \vert^{r}&\le 1+r\sum_{n\ge 1}\,\frac{1}{n^{1-r}}\,
\P( \vert {b}_{1} \vert>n)\\[2mm]
&\le 
 1+r\sum_{n\ge 1}\,\frac{1}{n^{1-r}}\,
\P\left( \sum^{l_{n}}_{j=1}\vert Q_{j} \vert>n\right)+
r\sum_{n\ge 1}\,\frac{1}{n^{1-r}}\,
\P( \vert t_{1} \vert>l_{n})
\\[2mm]
&\le 
 1+r\,\E\,\vert Q_{1} \vert^{\delta}\,
 \sum_{n\ge 1}\,\frac{l_{n}}{n^{1-r+\delta}}\,
 +r\,C\,
 \sum_{n\ge 1}\,\frac{1}{n^{1-r}l^{1/2}_{n}}\,.
\end{align*}

Putting here $l_{n}=[n^{4r}]$ we obtain that $\E\,\vert {b}_{1} \vert^{r}<\infty$ for any $r\in ]0,\delta/5[$.
Therefore, due to Proposition  \ref{Pr.sec:A.1}, we obtain that 
 for any bounded uniformly continuous function $f$
$$
\P\hbox{-}\lim_{N}\frac 1N\sum^{N}_{k=1}\,f(y_{k})=
\E f(\zeta),
$$
where 
\begin{equation}
\label{sec:Up.9}
\zeta={b}_{1}\,+\sum^\infty_{k=2}\,{b}_{k}\,
 \prod^{k-1}_{j=1}\, {a}_{j}
 \,.
\end{equation}
Now we show that 
$$
\P(\zeta<-x)>0
\quad\mbox{for any}\qquad 
x>0.
$$
Indeed, the random variable (\ref{sec:Up.9}) can be represented as
$$
\zeta={b}_{1}\,+{a}_{1}\,\zeta_{1},  \qquad    \zeta_{1}=\sum^\infty_{k=2}{b}_{k}\,
 \prod^{k-1}_{j=2}{a}_{j}.
$$

It is clear that the $\zeta_{1}$ is independent on ${b}_{1}$ and ${a}_{1}$. Note that on the set 
$\{\Delta \theta_1=1\}$ we have ${a}_{1}=M_{1}$ and ${b}_{1}=Q_{1}$. Therefore, for any $x>0$
$$
\P(\zeta<-x)\ge 
\P({b}_{1}\,+{a}_{1}\,\zeta_{1}<-x,\;\Delta \theta_1=1)\\[2mm]
= \P(Q_{1}\,+M_{1}\,\zeta<-x\,,\,w_{T_{1}}<0)\\[2mm]
$$
and we conclude as in the previous proposition.   \fdem

\section{Regularity of the ruin probability}
\label{sec:Reg}
\subsection{Integral representations}
The proof of smoothness of a function $H$ admitting an  integral representation  is based on a simple idea which merits to be explained.  

First, we recall the classical result on differentiability of the integral $H(u)=\int f(u,z)dz$ where $f(u,.)\in L^1$ for each $u$ from open subset $U\in \bbr$. {\it If 
$f(.,z)$ is differentiable on an open interval $]u_0-\e, u_0+\e[$ for almost all $z$ and,  on this interval, $|\partial f(.,z)/\partial u| \le g(z)$ (a.e.) where $g\in L^1$. Then $H$ is differentiable at $u_0$ and $H'(u_0)=\int \partial f(u_0,z)/\partial  u\, dz$}.    

Suppose that we are given 
a bounded measurable function $h(z)$ and a Gaussian random variable $\zeta \sim N(0,1)$.  Let 
$
H(u)=\E h(u+\zeta)=\int\,h(u+x)\varphi_{0,1} (x)dx.  $   
Then $H$  is differentiable  and even of the class $C^\infty$.  Of course, the above result cannot be applied directly. But using the change of variable we get the representation 
$$
H(u)=\int\,h(u+x)\varphi_{0,1} (x)dx= \int h(z)\varphi_{0,1} (z-u)dz. 
$$ 
Now the parameter $u$ appears only in the function $\varphi_{0,1}$, the integrand is differentiable in $u$ and we can apply the classical sufficient condition.   

 The issues here are: an integral representation, the smoothness of the density, and the integrability of its derivatives. 
   In the case  of survival probability $\Phi$  the integral representation is obtained from the strong Markov property. Unfortunately, the structure  of the representation 
is rather complicated,  the random variable standing for $\zeta$ is not a  Gaussian one, and its density is not given by a simple formula. Nevertheless, the idea of using a change of variable to move the 
parameter from the unknown function on which we  have only a limited  information (essentially, the boundedness and measurability)  
to the density  still does work. The main difficulty is to check the smoothness of the density and find appropriate bounds for its  derivatives. 
\medskip

 \begin{theorem}
\label{Th.sec:Reg.1}
Assume that the distribution function of $\xi_{1}$ has a  density $f$ differentiable on $\bbr_+$ and such that  $f'\in L^1(\bbr_{+})$.
 Then $\Phi (u)$ is two times continuously differentiable on  $]0,\infty[$.
 \end{theorem}
{\sl Proof.}  
We again consider the process
\begin{equation}
\label{sec:Reg.0-1}
Y^{u}_{t}=e^{\eta_t}\left (u
-R_t\right), 
\end{equation}
where  $R_t$ is defined in (\ref{defR}). 
Put 
\begin{equation}
\label{sec:Reg.0-0}
\theta^{u}:=\inf\{t\ge 0\colon \ Y^{u}_{t}\le 0\}. 
 \end{equation} 
 By virtue of  the strong Markov property of $X ^u$
\begin{equation}
\label{sec:Reg.0-1}
\Phi(u)=\E\Phi(X^u_{\theta^{u}\wedge T_1}).
 \end{equation} 
 
Note that the process $Y^u$ is strictly positive before the time $\theta^u$, zero at $\theta^u$, and strictly negative afterwards. Due to the independence of the Wiener process and the  instants of jumps, $\theta ^u\neq T_1$ a.s. Thus,  
$\{Y^u_{T_1}>0\}=\{\theta^u>T_1\}$ a.s.  Taking into account that $\Phi(0)=0$, we get that 
$$
\Phi(u)=\E I_{\{Y^u_{T_1}>0\}}\Phi(X^u_{T_1})=\E I_{\{Y^u_{T_1}>0\}}\Phi(Y^u_{T_1}+\xi_1)=\Phi_1(u)+\Phi_2(u), 
$$
where
$$
\Phi_1(u):=\alpha \int_0^2 \E G(Y_t^ u)e^{-\alpha t}dt,\qquad \ 
\Phi_2(u):=\alpha \int_2^\infty \E G(Y_t^ u)e^{-\alpha t}dt  
$$
with 
$$
G(y):=
\Chi_{\{y> 0\}}
\E\,\Phi\left(y+\xi_1\right)=\Chi_{\{y> 0\}}\int_0^\infty \Phi(y+x)dF(x).
$$
We analyze separately the smoothness of $\Phi_1$ and $\Phi_2$  using for this function 
appropriate integral representations. 

\subsection{Smoothness of $\Phi_2$}
We start with a simpler case of $\Phi_2$ and show that this function is infinitely differentiable without any assumptions on the distributions of $\xi_1$. 

From the representation 
$$
Y_t^u= e^{\eta_t-\eta_1}Y_1^u-c \int_1^t e^{\eta_t-\eta_s} ds, \qquad t\ge 1,  
$$
we obtain, using the independence of $Y_1^u$ and the process $(\eta_s-\eta_1)_{s\ge 1}$, that 
$$ 
\E (G(Y^u_t)|Y^u_1)=G(t,Y_1^u),
$$
where
$$
G(t,y) := \E G\left(e^{\eta_{t}-\eta_{1}}y-
c\int_{1}^{t}
e^{ \eta_{t} - \eta_{s}}d s\right). 
 $$ 
Substituting the expression for $Y^u_1$ given by (\ref{sec:Reg.0-1}) we have: 
$$
\Phi_2(u)=\E\int_2^{\infty}\E(G(Y^u_t)|Y^u_1)\alpha e^{-\alpha t}dt = 
\E H\left(e^{\kappa+\sigma W_{1}}(u-R_1) \right),
$$
where $H$ is a function taking values in $[0,1]$ and given by the formula 
$$
H(y):=\alpha\,\int_{2}^{\infty}G(t, y)e^{-\alpha t}d t.
 $$
 Taking into account that the conditional distribution of the process $(W_s)_{s\le 1}$ given $W_1=x$ is the same as of the Brownian bridge $B^x=(B_s^x)_{s\le 1}$ with  $B^x_s =W_s+s(x-W_1)$ we obtain the representation 
\beq
\label{rep3}
\Phi_2(u)=\int\E\, H(e^{\kappa+\sigma x}(u - \zeta^x))\varphi_{0,1}(x)dx 
\eeq
where
\begin{equation}
\label{zetax}
\zeta^x:=c\int_0^1e^{-(\kappa s +sx +\sigma  (W_s-sW_1))}ds.
\end{equation}
Lemma \ref{F2} below asserts that for every  $x$ the random variable $\zeta^x$ has a density $\rho (x,.)$ on $]0,\infty[$ and we easily obtain from 
(\ref{rep3}) by changing variable that 
\beq
\label{repr2}
\Phi_2(u)=\int_{0}^u\int H(e^{\kappa+\sigma x} z) \rho(x,u - z)\varphi_{0,1}(x)dxdz. 
\eeq

\begin{lemma}
\label{F2}
The random variable $\zeta^x$ has a density  $\rho(x,.)\in C^\infty$ such that  for any $n\ge 1$  
\beq
\label{bderivs}
\sup_{y\ge 0}\,
\left\vert \frac {\partial^n}{\partial y^n} \rho(x,y) \right\vert  \le C_n e^{C_n |x|}
\eeq
with some constant $C_n$ and  $(\partial^n /\partial y^n)\rho(x,0)=0$. 
\end{lemma}
\noindent
{\sl Proof.} We obtain the result using again the integral representation. 
Let us introduce the random process
\bean
D_s&:=&\big((W_s-2sW_{1/2})+s(W_{1/2}-W_1)\big)I_{\{s\le 1/2\}}\\
&& +\big((1-s)(W_s-W_{1/2}))-s(W_1-W_s)\big)I_{\{s>1/2 \}}, 
\eean
and the piecewise linear function 
$$
\gamma_s:=sI_{\{s\le 1/2\}}+(1-s)I_{\{s>1/2 \}},  \quad s\in [0,1].  
$$
The following identity is obvious: 
$$
W_s-sW_1=D_s+  \gamma_sW_{1/2}. 
$$
Since $D_s$ and $W_{1/2}$ are independent random variables, for any bounded Borel function $g$ we have: 
$$
\E\,g(\zeta^x)=\E\int g\Big(c\int_0^1e^{-(\kappa s +sx + \sigma  D_s+ \sigma  \gamma_sv)}ds\Big)\varphi_{0,1/2}(v)dv
$$
Let $v(x,.)$ be the inverse of the continuous strictly decreasing function 
$$
y\mapsto c\int_0^1e^{-(\kappa s +sx +\sigma D_s+\sigma  \gamma_sy)}ds
$$
depending on the parameter $x$ (and also  on $\omega$ omitted as usual).  Note that $v(x,0+)=\infty$ and $v(x,\infty)=0$. 
After the change of variable, we obtain, using the notation 
$$
K(x,z):=c\sigma\int_0^1\gamma_se^{-(\kappa s +sx +\sigma D_s+ \sigma \gamma_s z)}ds, 
$$
 that 
$$
\E g(\zeta^x)=\int^{\infty}_0 g(y)\rho(x,y) dy,  
$$
where  
\beq
\rho(x,.) := \E \frac {\varphi_{0,1/2}(v(x,.))}{K(x,v(x,.))}.   
\eeq
Thus, $\rho(x,.)$ is the density of distribution of the random variable $\zeta^x$. It remains to check that it is infinitely differentiable and find appropriate bounds for its derivatives.  
 
Put 
$$
Q^{(0)}(x,z):=\frac {\varphi_{0,1/2}(z)}{K(x,z)}, \qquad 
 Q^{(n)}(x,z):=-\frac{Q^{(n-1)}_z(x,z)}{K(x,z)}, \quad n \ge 1.  
$$
Then 
$$
\frac{\partial }{\partial y}\,Q^{(0)}(x,v(x,y))=
Q^{(0)}_z(x,v(x,y))v_y(x,y)=
-
\frac{Q^{(0)}_z(x,v(x,y))}{K(x,v(x,y))}
=Q^{(1)}(x,v(x,y))  
$$
and, similarly, 
$$
\frac{\partial^n }{\partial y^n}\,Q^{(0)}(x,v(x,y))=
-
\frac{Q^{(n-1)}_z(x,v(x,y))}{K(x,v(x,y))}
=Q^{(n)}(x,v(x,y)).
$$
It is easily seen that 
\beq
\label{Qn}
Q^{(n)}(x,z)=\frac{\varphi_{0,1/2}(z)}{K^{n+1}(x,z)}\sum_{k=0}^nP_k(z)R_{n-k}(x,z),  
\eeq
where $P_k(z)$ is a polynomial of order  $k$ and $R_{n-k}(x,z)$ is a linear combination of products 
of derivatives of  $K(x,z)$ in variable $z$. 
Note that for any $x,y\in [0,1]$ we have the bounds
$$
{-|x| -\sigma |z| - 3\sigma W^*_1}\le  {\kappa s +sx +\sigma D_s+ \sigma  \gamma_s z}\le {\kappa  +|x| +\sigma |z| +3\sigma W^*_1}
$$
where $ W^{*}_1=\sup_{s\le 1}\,|W_{s}|$.
It follows that   there  exists a constant 
$C_n>0$  such that
\beq
\label{cricialbound}
\vert Q^{(n)}(x,z)\vert
\le C_n e^{C_n |x|}(1+z^n)e^{-z^2} e^{C_n\,W^*_1}. 
\eeq
Since  $\E\,e^{C_n\,W^*_1}<\infty$,  for each $x,y$ and $n$ the derivative 
$({\partial^n }/{\partial y^n})\,Q^{(0)}(x,v(x,y))$ admits a $\P$-integrable bound.  
Thus,  we can differentiate under the sign of expectation and obtain that 
$$
\frac{\partial^n }{\partial y^n}\rho (x,y)=\E \frac{\partial^n }{\partial y^n}\frac {\varphi_{0,1/2}(v(x,y))}{K(x,v(x,y))}
=
\E \frac{\partial^n }{\partial y^n}\,Q^{(0)}(x,v(x,y))=\E Q^{(n)}(x,v(x,y)).    
$$ 
Moreover, the bound  (\ref{cricialbound})
ensures that,  for some constant $\tilde C_n$,
$$
\sup_{y\ge 0}
\E\left \vert \frac{\partial^n }{\partial y^n}\,Q^{(0)}(x,v(x,y)) \right\vert \,
\le\,
\E\,\sup_{z\in\bbr}
\vert Q^{(n)}(x,z)\vert
\le \tilde C_n e^{C_n\vert x\vert}
$$
and the bound (\ref{bderivs}) holds. 

Since $v(x,0+)=\infty$, the bound (\ref{cricialbound}) implies that  $(\partial^n /\partial y^n)\rho(x,0)=0$. 
\fdem

\begin{remark}
It is worth to trace in these arguments  the dependence of the constant $\tilde C^n$ on the  parameters $c$ and $\sigma$ when they are approaching zero. From the formula (\ref{Qn})
it is clear that $\tilde C_n$ should be proportional to $(c\sigma)^{-n}$. 
\end{remark}

\begin{proposition}
\label{Phi2.1}
The function $\Phi_{2}(u)$ belongs to $C^{\infty}(]0,\infty[)$.
\end{proposition}
{\sl Proof.}  Putting 
$$
\wt{H}(u,z):=\int H(e^{\kappa+\sigma x} z) \rho(x,u - z)\varphi_{0,1}(x)
\,dx,
$$
we rewrite the formula (\ref{repr2}) as 
\begin{equation}
\label{phi2.2}
\Phi_2(u)= \int^{u}_{0}\,\wt{H}(u,z)\,
dz. 
\end{equation}
Clearly, the function $\wt{H}$ is continuous  on $]0,\infty[\times ]0,\infty[$.  Using Lemma \ref{F2} we obtain that 
$$
\frac{\partial}{\partial u}\,\wt{H}(u,z)\,
=
\int  H(e^{\kappa+\sigma x} z) \frac{\partial}{\partial u}\, \rho(x,u - z)
\varphi_{0,1}(x)
\,dx
$$
and
$$
\sup_{u,z}\,
\left\vert
\frac{\partial}{\partial u}\,\wt{H}(u,z)\,
\right\vert
<\,\infty.
$$
By induction, for every  $n\ge 1$, 
$$
\frac{\partial^{n}}{\partial u^{n}}\,\wt{H}(u,z)
=
\int  H(e^{\kappa+\sigma x} z)
\frac{\partial^{n}}{\partial u^{n}}\, \rho(x,u - z)
\varphi_{0,1}(x)
\,dx
$$
and
$$
\sup_{u,z}\,
\left\vert
\frac{\partial^{n}}{\partial u^{n}}\,\wt{H}(u,z)\,
\right\vert
\,<\,\infty.
$$
By virtue of  Lemma \ref{F2} 
$\rho(x,0)=0$, i.e. $\wt{H}(u,u)\,=0$. 
So,
$$
ù\label{phi2.3}
\frac {d}{du}\Phi_2(u)
=\, \wt{H}(u,u)\,+
 \, \int^{u}_{0}\,
 \frac{\partial}{\partial u}\,
 \wt{H}(u,z)\,
dz
=
 \, \int^{u}_{0}\,
 \frac{\partial}{\partial u}\,
 \wt{H}(u,z)\,dz. 
$$
In the same way we  check that
$$
\frac {d^n}{du^n}\Phi_2(u)
=
\, \int^{u}_{0}\,
 \frac{\partial^{n}}{\partial u^{n}}
 \wt{H}(u,z)\,
dz. 
$$
for any $n\ge 1$. 
\fdem

\subsection{Smoothness of $\Phi_1$}
Arguing in the same spirit as in the previous subsection but taking this time the conditional expectation with respect to $W_t$  we obtain that 
$$
\E \Chi_{\{R_{t}<u\}}\,h(e^{\kappa {t}+\sigma W_t}(u-R_{t}))=
\frac{1}{\sqrt{t}}\E
\int
\Chi_{\{\zeta^{t,x}<u\}}
h(u,t,x)\varphi_{0,1}\left(\frac {x}{\sqrt{t}}\right)dx
$$
where we use the abbreviations 
$$
h(u,t,x):=h(e^{\kappa t+\sigma x}
(u-\zeta^{t,x})),\qquad 
h(y)=\E\,\Phi\left(y+\xi_1\right)
$$
and
$$
\zeta^{t,x}:=c\int^{t}_{0}e^{- (s x/t+\kappa s+
\sigma (W_{s}-(s/t)W_{t})
}\,d s. 
$$

It is easily seen that the random variable $\zeta^{t,x}$ has infinitely differentiable density (the same as of 
$\zeta^x$ defined in (\ref{zetax}) but with the parameters $ct$, $\kappa t$, and $\sigma t^{1/2}$). 
Unfortunately,  derivatives of this density have non-integrable singularities as $t$ tends to zero (see Remark at the end of previous  subsection).  By this reason  we cannot use the strategy of proof used for $\Phi_2$. Nevertheless, the hypothesis on the distribution of $\xi_1$ allows us to establish the claimed result. 

Note that the function  $x\to \zeta ^{t,x}$ is strictly decreasing and maps $\bbr$ onto $\bbr_+$. Let denote 
 $z(t,.)$ its inverse.   The derivative of the latter is given by the formula
$$
z_x(t,x)=-\frac{t}{L(t,z(t,x))},
$$
where
\begin{equation}
\label{sec:Reg.7}
L(t,z)=c \int^{t}_{0}se^{- (s z/t+\kappa s+
\sigma (W_{s}-(s/t)W_{t})}\,ds.
\end{equation}
Changing the variable 
we obtain that
$$
\int
\Chi_{\{\zeta^{t,z}<u\}}
{h}(u,t,z)
\varphi_{0,1}\left(\frac{z}{\sqrt{t}}\right)d z
=t
\int^{u}_{0}
{h}(u,t,z(t,x))
\,
D(t,z(t,x))\,
d x,
$$
where 
$$
D(t,z)=\frac{\varphi_{0,1}\left({z}/\sqrt{t}\right)}{L(t,z)}.
$$
Summarizing, we get that 
$$
\Phi_{1}(u)=\alpha 
\,\E\,\int^{2}_{0}\,
\sqrt{t}\,
H(t,u)\,e^{-\alpha t}\,d t,
$$
where
\begin{equation}
\label{sec:Reg.9}
H(t,u):=
\int^{u}_{0}
{h}(u,t,z(t,x))
\,
D(t,z(t,x))\,
d x.
\end{equation}


\begin{proposition}
 \label{Pr.Reg.3}
Under the conditions of Theorem~\ref{Th.sec:Reg.1} the function $H(t,u)$ defined in 
 (\ref{sec:Reg.9})
 for any fixed $u_{0}>0$ 
satisfies the following inequality
  \begin{equation}
\label{sec:Reg.9-10}
\sup_{t\in]0, 2]}\,
\E\,
\sup_{u\ge u_{0}}\,
\left(\vert
H_u(t,u)\vert
+
\vert
H_{uu}(t,u)
\vert
\right)
<\infty.
\end{equation}
  \end{proposition}
{\sl Proof.}  In virtue of the hypothesis the function  $h(y)=\E \Phi(y+\xi_{1})$ is differentiable.  Differentiating (\ref{sec:Reg.9}) we get that 
\begin{equation}
\label{sec:Reg.10-0}
H_u(t,u)=h(0)\,D(t,z(t,u))\,+
\int^{u}_{0}h_u(u,z(t,x))\,
D(t,z(t,x))\,
d x,
\end{equation}
where $h(0)=\E \Phi(\xi_1)$ and 
$$
h_u(u,t,z(t,x))
=h'\big(e^{\kappa t+\sigma z(t,x)}(u-x)\big)
e^{\kappa t+\sigma z(t,x)}.
$$
Note that  
$$
\frac{\partial }{\partial x}h(u,t,z(t,x))=h_u(u,t,z(t,x))\big[\sigma z_x(t,x)(u-x)-1\big].
$$
Therefore, 
$$
\int^{u}_{0}h_u(u,z(t,x))
D(t,z(t,x))\,d x
=
-\int^{u}_{0}\,
\frac{\varphi_{0,1}\left({z(t,x)}/{\sqrt{t}}\right)}{\sigma t(u-x)+L(t,z(t,x))}
\,d_{x}\, h(u,t,z(t,x)).
$$
Integrating by parts and taking into account that $z(t,0+)=\infty$ we get  that 
\begin{equation}
\label{sec:Reg.10-01}
H_u(t,u)=
\int^{u}_{0}
h(u,t,z(t,x))\,
\Theta (u,t,z(t,x))\,\varphi_{0,1}\left({z(t,x)}/{\sqrt{t}}\right)\,
d x, 
\end{equation}
where 
\begin{align*}
\Theta(u,t,z)&
=\frac{z}{L(t,z)(\sigma t\left(u-\zeta^{t,z})\right)+L(t,z))}
-
\frac{t\sigma L(t,z)+tL_z(t,z)}{L(t,z)\left(\sigma t\left(u-\zeta^{t,z}\right)+L(t,z)\right)^{2}}. 
\end{align*}

Inspecting the formula (\ref{sec:Reg.7}) defining $L(t,z)$ we conclude that there exist positive constants $C_0$ ("small") and $C_1$ ("large") such that 
$$
\max_{t\in ]0,2]}\,
\left( 
L(t,z)+
t|L_z(t,z)|
\right)
\le C_1
e^{C_1(|z|+W^{*}_{t})}
$$
and for any $t\in ]0,2]$
$$
|L(t,z)|\ge C_{0}\,t^{2}\,
e^{- C_{1}(|z|+W^{*}_{t})},
$$
where $W^{*}_{t}:=\max_{v\le t}|W_{v}|$.
Taking this into account we obtain that
for some   $C>0$
\begin{equation}
\label{sec:Reg.11}
\left\vert
\Theta(u,t,z)
\right\vert
\le C\,t^{-6}\,
e^{C(|z|+W^{*}_{t})}
\,.
\end{equation}
Using a generic notation $C$ for a constant (which may vary even within a single formula) we obtain, for any $t\in]0, 2]$ 
\begin{align*}
\left\vert
H_u(t,u)
\right\vert
&\le
Ct^{-7}e^{CW^{*}_{t}}\int^{\infty}_{z(t,u)}\,e^{C|z|-\frac{z^{2}}{2t}}L(t,z)\, d z
\\[4mm] 
&\le Ct^{-7}
e^{C\,W^{*}_{2}}\,e^{-\frac{z^{2}(t,u)}{4t}}\,
\int^{\infty}_{0}\,e^{C|z|-\frac{z^{2}}{4}}\,d z. 
\end{align*}
So, we have the bound 
$$
\left\vert
H_u(t,u)
\right\vert
\le C t^{-7}
e^{CW^{*}_{2}}\,e^{-\frac{z^{2}(t,u)}{4t}}.
$$
For any $u\ge u_{0}>0$ and $t\in ]0,2]$ 
$$
u_{0}\le u=\zeta^{t,z(t,u)}\le c t\, e^{2|\kappa|+\sigma |z(t,u)|+ \sigma W^{*}_{t}},
$$
i.e.
$$
|z(t,u)|\ge \frac{1}{\sigma}\,\ln\,\frac{e^{-2|\kappa|}u_{0}}{tc}-W^{*}_{t}.
$$
Put
$$
t_{0}:=\,\min\left(\frac{e^{-2|\kappa|-3\sigma} u_{0}}{c},\,2\right), 
\qquad
\Gamma:=\{W^{*}_{t}\le\,t^{1/4}\}.
$$
Thus, for $t\in ]0,t_{0}]$ on the set $\Gamma$ we have the inequality 
$$
|z(t,u)|\ge 1.
$$
Taking into account that $\E\,e^{aW^{*}_{t}}<\infty$ for any $a$ and $t>0$,  
we obtain that
$$
\E\,
\max_{t\in [t_{0},2]}\,\sup_{u\ge 0}\,\left\vert
H_u(t,u)
\right\vert\,
\le C t_{0}^{-7}\,\E\,
e^{CW^{*}_{2}}
<\infty.
$$
For $t\in ]0,t_{0}]$  we have
\begin{align*}
\E\,\max_{u\ge u_{0}}\,\left\vert
H_u(t,u)
\right\vert\,&
\le 
C t^{-7} 
\left(
e^{-\frac{1}{4t}}
+
\E e^{CW^{*}_{t}}\,\Chi_{\Gamma^{c}}
\right)\\[4mm]
&\le 
C t^{-7} 
\left(
e^{-\frac{1}{4t}}
+
\,\sqrt{\E\,e^{CW^{*}_{2}}} 
\sqrt{\P\left(W^{*}_{t}\ge t^{1/4}\right)}
\right).
\end{align*}
By the Chebyshev inequality we have: 
$$
\P\left(W^{*}_{t}\ge t^{1/4}\right)\le e^{-t^{-1/4}}\,\E\,e^{\frac{W^{*}_{t}}{\sqrt{t}}}=e^{-t^{-1/4}}\,\E\,e^{W^{*}_{1}}, 
$$
that is 
$$
\sup_{t\ge 0} e^{t^{-1/4}}\P\left(W^{*}_{t}\ge t^{1/4}\right)<\infty.
$$
This implies that
$$
\max_{t\in]0, t_{0}]}
\E\max_{u\ge u_{0}}\left\vert
H_u(t,u)
\right\vert
<\infty.
$$
Therefore,
$$
\max_{t\in]0, 2]}
\E\max_{u\ge u_{0}}\,\left\vert
H_u(t,u)
\right\vert
\,<\,\infty.
$$

Differentiating  (\ref{sec:Reg.10-01}) we find that
$$
H_{uu}(t,u)=
\,h(0)\,
\Theta (u,t,z(t,u))\,\varphi\left(\frac{z(t,u)}{\sqrt{t}}\right)\,
+
\int^{u}_{0}
\Upsilon(u,t,z(t,x))\,
\varphi\left(\frac{z(t,x)}{\sqrt{t}}\right)\,
d x,
$$
where
$$
\Upsilon (u,t,z)=
{h}_{u}(u,t,z)
\Theta(u,t,z)
+
{h}(u,t,z)
\Theta_u(u,t,z).
$$
By assumption, the distribution function $F$ has the density $f$ whose derivative $f'$ is a continuous function on $\bbr_+$ integrable with respect to the Lebesgue measure.   
Changing variable we get that  
\begin{align*}
{h'}(y)&=\frac {d}{dy}\int_0^{\infty}\Phi(y+x)f(x)dx=
\frac {d}{dy}\int_y^{\infty}\Phi(x)f(x-y)dx\\
=&-\Phi(y)\,f(0)
-\,
\int^{\infty}_{y}\,\Phi(z){f'}(z-y)\,d z,
\end{align*}
i.e.
$$
\sup_{x\ge 0}\,|{h'}(x)|\,<\,\infty\,.
$$
Similarly to (\ref{sec:Reg.11})  we obtain that
$$
\left\vert
\Upsilon_{t}(u,z)
\right\vert
\,\le\,C\,t^{-8}\,
e^{C(|z|+W^{*}_{t})}.
$$
Therefore, 
$$
\left\vert
H_{uu}(t,u)
\right\vert\,
\le\,C t^{-9}\,
e^{C\,W^{*}_{t}}\,e^{-\frac{z^{2}(t,u)}{4t}}, 
$$
implying that  
$$
\max_{t\in]0, 2]}\,
\E\,\max_{u\ge u_{0}}\,\left\vert
H_{uu}(t,u)
\right\vert
<\infty.
$$
Proposition~\ref{Pr.Reg.3} is proven.
\fdem

\subsection{Integro-differential equation  for the survival probability}
\label{sec:Eq}
\begin{proposition} 
\label{Pr.Eq.1}
Suppose that $\Phi\in C^2$. Then 
$\Phi$ satisfies the equation (\ref{int-diff}).
  \end{proposition}
{\sl Proof.} 
For $h>0$ and $\epsilon>0$ assumed to be small enough to ensure that $u\in ]\epsilon,\epsilon^{-1}[$ we put 
$$
\tau^{\epsilon}_h:=\inf\big\{t\ge 0\colon\  Y^{u}_{t}
\notin[\epsilon\,,\epsilon^{-1}]\big\}\wedge h. 
$$
 Using the Ito formula and taking into account that on the interval on $[0,T_1[$ the process $X^u$ coincides with $Y^u$   we obtain the representation 
 \bean
 \Phi(X^u_{\tau^{\epsilon}_h\wedge T_1})
 &=&
 \Phi(u)+\sigma \int_{0}^{\tau^{\epsilon}_h\wedge T_1}\Phi'(Y^u_{s})\,dW_s\\
 &&+ \int_{0}^{\tau^{\epsilon}_h\wedge T_1}\Big(\frac 12{\sigma^2}
  (Y^u_{s})^2\Phi''(Y^u_{s})+(aY^u_{s}-c)\Phi'(Y^u_{s})\Big)ds\\
 &&+
(\Phi(Y^u_{T_1}+\xi_1)-\Phi(Y^u_{T_1}))I_{\{T_1\le {\tau^{\epsilon}_h}\}}.
 \eean
Due to the strong Markov property $\Phi(u)=\E\,\Phi(X^u_{\tau^{\epsilon}_h\wedge T_1})$. For every $\epsilon>0$ the integrands above are 
bounded by  constants and, hence, the expectation of the stochastic integral is zero.  Moreover, $\tau^\epsilon_h\wedge T_1=h$ when $h$ is sufficiently small  
(the threshold below which we have this equality, of course, depends on $\omega$). 

It follows that, independently of $\epsilon$,  
$$
\frac 1h \E\int_{0}^{\tau^{\epsilon}_h\wedge T_1}\Big(\frac 12{\sigma^2}
  (Y^u_{s})^2\Phi''(Y^u_{s})+aY^u_{s}\Phi'(Y^u_{s})-c\Big)ds\to \frac 12{\sigma^2}
  u^2\Phi''(u)+(au-c)\Phi'(u). 
$$
Finally,  
$$
\frac 1h \E(\Phi(Y^u_{T_1}+\xi_1)-\Phi(Y^u_{T_1}))I_{\{T_1\le {\tau^{\epsilon}_h}\}}=\alpha \E \frac 1h \int_{0}^{\tau^{\epsilon}_h} 
(\Phi(Y^u_{t}+\xi_1)-\Phi(Y^u_{t}))e^{-\alpha t}\,dt. 
$$
The right-hand converges to $\alpha(\E\Phi(u+\xi_1)-\Phi(u))$ as $h\to 0$ in virtue of the Lebesgue theorem on dominated convergence. 
It follows that $\Phi$ satisfies the equation (\ref{int-diff}). \fdem

\section{Proof of Theorem~\ref{main}}

 Assume that the claims are exponentially distributed, i.e.
 $F(x)=1-e^{-x/\mu}$. Similarly to the classical case, this assumption
allows us to obtain for the ruin probability an ordinary differential
 equation (but of a higher order). Indeed, now
 the integro-differential  equation (\ref{int-diff}) has the form 
\begin{equation}
\label {int-diff1} \frac 12 \sigma^2u^2 \Phi''(u)+ (au-c)\Phi'(u) -
\alpha\Phi(u)+\frac {\alpha}{\mu} \int_0^\infty\Phi(u+y)e^{-y/\mu}dy=0.
\end{equation}
 Notice that
 $$
  \frac{d}{du}\int_0^\infty\Phi(u+y)e^{-y/\mu}dy=-\Phi(u)
  +
  \frac 1\mu \int_0^\infty\Phi(u+y)e^{-y/\mu}dy.
 $$
 Differentiating  (\ref{int-diff1}) and adding to it the obtained identity multiplied by $\mu$ 
 we exclude the integral term and arrive at a third order ordinary differential equation.
 It does not contain the function
 itself and, therefore, is reduced to a second order differential equation  for the function $G=\Phi'$ which can be easily transformed to the form 
\begin{equation}
\label{diff}
G^{\prime\prime}- p(u)G^{\prime}+ p_{0}(u)G=0,
\end{equation}
where
\begin{align*}
p(u)&:=\frac{1}{\mu}-
2\left (1+\frac{a}{\sigma^2}\right)
\frac{1}{u}+\frac {2c}{\sigma^2}
\frac{1}{u^2},\\[4mm]
 p_0(u)&:=-
\frac{2a}{\mu\sigma^2}\frac{1}{u}+\left(a-\alpha +c/\mu\right)
 \frac {2}{\sigma^2}\frac 1{u^2}.
\end{align*}
The substitution $G(u)=R(u)Z(u)$ with
$$
R(u):=\exp\left\{\frac{1}{2} \int_{1}^u p(s)\,d s\right\}
$$
eliminates the first derivative and leads to the equation
\begin{equation}
\label{diff-1}
Z''- q(u) Z=0\,,
\end{equation}
where
$$
q(u):=\frac{1}{4\mu^{2}}+\left(\frac a{\sigma^2}-1\right)\frac{1}{u}+\sum _{i=2}^4
A_{i}\frac{1}{u^i}
$$
with certain constants $A_{i}$ which are of no importance in our asymptotic analysis. It is easy to check that
$$
\int^{\infty}_{x_{0}}\,\left( 
\frac{|q^{\prime\prime}(x)|}{q^{3/2}(x)}
+
\frac{|q^{\prime}(x)|^{2}}{q^{5/2}(x)}
\right)\,d x\,<\,\infty,
$$
where $x_{0}=\sup\{x\ge 1\colon\ q(x)\le 0\}+1$.

According
to \cite{Fed}, pp. 54-55, the equation \eqref{diff-1}
 has two  fundamental solutions $Z_{+}$
and $Z_{-}$ 
$$
Z_{\pm}(x)=\sqrt{2\mu}\,\exp\left\{\pm\int^{x}_{x_{0}}\,\sqrt{q(z)}d z\right\}\left(1+o(1)\right), \qquad x\to\infty,
$$
i.e.
$$
Z_{\pm}\sim \exp\left\{\pm\left(\frac{x}{2\mu}\,+\,\frac{a-\sigma^{2}}{\sigma^{2}}\,\ln x\right)\right\}.
$$
Since 
$$
R(x)\sim \exp\left\{\frac{x}{2\mu}-\frac{a+\sigma^{2}}{\sigma^{2}}\,\ln x \right\},
$$
we obtain that \eqref{diff} admits, as
solutions,  functions with the following asymptotics:
$$
G_{+}(x)\sim x^{-2}e^{\frac{1}{\mu}x},
\qquad
G_{-}(x)\sim x^{-2a/\sigma^2} .
$$
The differential equation of the third order for  $\Phi$ has three
solutions: $ \Phi_{0}(u)=1$,
$$
\Phi_{+}(u)=\int^{u}_{x_{0}}\,G_{+}(x)\,d x
\qquad
 \Phi_{-}(u)=-\int^{+\infty}_{u}\,G_{-}(x)\,d x\,.
$$
 The ruin probability
$\Psi:=1-\Phi$ is the linear combination of these functions, i.e.
$$
\Psi (u)=C_{0}+C_{1}\Phi_{+}(u)+C_{2}\Phi_{-}(u).
$$
Since $\Phi_{+}(u)\to\infty$ as $u\to\infty$, we obtain immediately that $C_{1}=0$. 
For the case $\beta>0$ we know from Proposition \ref{main-0}
that $\Psi (\infty)=0$. Thus, $C_0=0$ and 
$$
\Psi (u)= C_{2}\,\int_{u}^\infty\, x^{-2a/\sigma^2}\,\left(1+\delta(x)\right)\,d x,
$$
where $\delta(x)\to 0$ as $x\to\infty$ and $C_2>0$. 
The  integral decreases at infinity as the power function
$u^{-\beta}/\beta$ and we obtain  Theorem
\ref{main}. \fdem

\renewcommand{\theequation}{A.\arabic{equation}}
\renewcommand{\thetheorem}{A.\arabic{theorem}}
\renewcommand{\thesubsection}{A.\arabic{subsection}}
\section{Appendix: ergodic theorem for an autoregression  with  random coefficients}\label{sec:A}
\setcounter{equation}{0}
\setcounter{theorem}{0}


Let $(a_{n}\,,b_{n})_{n\ge 1}$ be an  i.i.d.
sequence of random variable in  $\bbr^2$ and let $x_0$ be an arbitrary constant.  Define  the sequence $(x_n)$ recursively by the formula  
$$
x_{n}\,=\,a_{n}\,x_{n-1}\,+\,b_{n}, \quad n\ge 1.
$$

\begin{proposition}\label{Pr.sec:A.1}
Assume that there exists $\delta \in ]0,1]$ such that 
$$
\E\,|a_{1}|^\delta<1,
\qquad
\qquad 
\E\,|b_{1}|^\delta<\infty. 
$$
Then
 for any bounded uniformly continuous function $f$
$$
\P\hbox{-}\lim_{N}\frac1N\sum^N_{n=1}\,f(x_n)=\
\E\,f(\zeta),
$$
where 
$$
\zeta=\sum^\infty_{k=1}\,b_{k}\,
 \prod^{k-1}_{j=1}\, a_{j}
 \quad\mbox{with}\quad
 \prod^{0}_{j=1}a_{j}=1.
 $$
\end{proposition}

The proof is given in \cite{PZ}.

\acknowledgement{The research is funded by the grant of the Government of Russian Federation \\ $n^\circ$14.А12.31.0007. 
The second author is partially
supported by the department of  Mathematics and Mechanics of National Research Tomsk State University
}.

\end{document}